\documentclass[12pt,a4paper]{article}
\usepackage{amssymb,amscd,amsmath,amsthm}
\usepackage{graphicx}
\usepackage{psfrag}
\theoremstyle{plain}
\newtheorem{thm}{Theorem}[section]
\newtheorem{prop}[thm]{Proposition}
\newtheorem{lem}[thm]{Lemma}
\newtheorem{cor}[thm]{Corollary}
\theoremstyle{definition}
\newtheorem{rem}[thm]{Remark}
\newtheorem{dfn}[thm]{Definition}
\newcommand{\claim}{{\bf Claim.}\quad}
\newcommand{\stepone}{{\bf Step 1.}\quad} 
\newcommand{\steptwo}{{\bf Step 2.}\quad} 
\newcommand{\stepthree}{{\bf Step 3.}\quad} 
\newcommand{\stepfour}{{\bf Step 4.}\quad} 
\newcommand{\stepfive}{{\bf Step 5.}\quad} 
 
\newcommand{\nandc}{{\bf Notations\ and\ Convention.}\quad} 
\newcommand{\BCC}{\mathbb{C}}
\newcommand{\BPP}{\mathbb{P}}
\newcommand{\BQQ}{\mathbb{Q}}
\newcommand{\BRR}{\mathbb{R}}
\newcommand{\BZZ}{\mathbb{Z}}

\newcommand{\MFD}{\mathfrak{D}}
\newcommand{\MCO}{\mathcal{O}}
\newcommand{\vp}{\varphi}
\newcommand{\Aut}{\operatorname{Aut}}
\newcommand{\algdim}{\operatorname{alg.dim}}

\newcommand{\Hom}{\operatorname{Hom}}
\newcommand{\id}{\operatorname{id}}
\newcommand{\im}{\operatorname{im}}
\newcommand{\mymod}{\hspace{-0.5mm}\operatorname{mod}}
\newcommand{\ord}{\operatorname{ord}}

\newcommand{\rank}{\operatorname{rank}}

\newcommand{\Spec}{\operatorname{Spec}}

\begin{document}
\title{On the number of Enriques quotients of a K3 surface\thanks{Communicated
 by S. Mukai. Received April 11, 2006, Revised May 17, 2006.
 2000 Mathematics Subject Classification: 14J28} }
\author{Hisanori Ohashi
\thanks{Research Institute for Mathematical Sciences,
 Kyoto University, Kyoto 606-8502, Japan. e-mail:pioggia@kurims.kyoto-u.ac.jp} } 
\maketitle
\date{}
\setcounter{section}{-1}
\section{Introduction}
\label{Intro}
 \quad A {\it $K3$ surface} $X$ is a compact complex surface 
 with $K_X\sim 0$ and $H^1(X,\MCO_X)=0$. An {\it Enriques surface} 
 is a compact complex surface with $H^1(Y,\MCO_Y)=H^2(Y,\MCO_Y)=0$
 and $2K_Y\sim 0$. The universal covering of an Enriques surface is 
 a $K3$ surface. Conversely every quotient of a $K3$ surface 
 by a free involution is an Enriques surface. Here a {\it free involution}
 is an automorphism of order $2$ without any fixed points.

 The moduli space of Enriques surfaces is constructed using the periods of 
 their covering $K3$ surfaces. Precisely speaking, an Enriques surface 
 determines a lattice-polarized $K3$ surface and vice versa, so that 
 the moduli space of Enriques surfaces can be described by the moduli 
 space of lattice-polarized $K3$ surfaces. We note that 
 even if we do not fix any polarization on Enriques surfaces, their 
 covering $K3$ surfaces automatically have a lattice-polarization.
 Then, what happens if we drop the lattice-polarization 
 of the covering $K3$ surface? 

 We will call two Enriques 
 quotients of a $K3$ surface {\it distinct} if they are not isomorphic to 
 each other as varieties. In his paper \cite{Kondo}, Kondo discovered a $K3$ surface
 with two distinct Enriques quotients. He computed the 
 automorphism groups of the two quotients. Since then, as far as 
 the author knows, no other examples have been found.

 In this paper we investigate this phenomenon.
 We show that $K3$ surfaces with more than one distinct Enriques
 quotients have $9$-dimensional components 
 (neither irreducible nor closed) in the period domain. Moreover we compute the 
 exact number of distinct Enriques quotients at a very general point
 of each component (Proposition \ref{div}).
 This generalizes Kondo's example in an arithmetic way and results
 in the following unboundedness theorem.
\begin{thm}\label{unbound}
 For any nonnegative integer $l$, there exists a $K3$ surface $X$ with 
 exactly $2^{l+10}$ distinct Enriques quotients. In particular, there does not 
 exist a universal bound for the number of distinct Enriques quotients of 
 a $K3$ surface.
\end{thm}
 \quad We also generalize Kondo's example in a geometric way.
 Its construction is due to Mukai \cite{Mukai}. We introduce his construction and show that 
 a generic Kummer surface $X$ of product type (see Section $4$) has exactly 
 $15$ distinct Enriques quotients, which can be constructed from classical Lieberman's
 involution and Kondo-Mukai's involution.
\begin{thm}
 \quad $X$ has exactly 15 distinct Enriques quotients which 
 are naturally in one-to-one correspondence with nonzero elements
 of the discriminant group of $NS(X)$.
\end{thm}
 \quad From the theoretical point of view, we first show the following finiteness theorem
 on the automorphism group of a $K3$ surface $X$ using a theorem of Borel.
\begin{thm}\label{FP}
 \quad In $\Aut(X)$ there are only finitely many conjugacy classes
 of finite subgroups.
\end{thm}
 This theorem concerns us because it bounds the number of distinct Enriques quotients 
 for any $K3$ surfaces.
\begin{cor}\label{fEq}
 Every $K3$ surface $X$ has only finitely many distinct Enriques quotients.
\end{cor} 
 The usage of the theorem of Borel is suggested by an anonimous referee.
 We remark that Corollary \ref{fEq} follows also from our counting method
 described in Section \ref{number}. There, 
 to count the exact number of distinct Enriques quotients, we consider more 
 directly the embeddings of the Enriques lattice $U(2)\oplus E_8(2)$
 into N\'{e}ron-Severi lattices.
 \vspace{5mm} \\
 $\nandc$\\
 \quad Our main tool is the theory of lattices and their discriminant forms.
 Here we collect some basic definitions about them. 
 See \cite{Nikulin1} for the detailed exposition.\\
 \quad A {\itshape lattice} $L$ is a free $\BZZ$-module of 
 finite rank equipped with a $\BZZ$-valued symmetric bilinear form.
 $L$ is said to be {\itshape even} if for all $l\in L, l^2\in 2\BZZ$.
 In this paper we treat only even lattices, so that we sometimes
 omit mentioning the evenness.
 For a lattice $L$, there is a natural homomorphism $c:L\rightarrow 
 L^*=\Hom (L,\BZZ)$ defined by $l\mapsto (l,\cdot)$. $L$ is said to be 
 {\itshape nondegenerate} if $c$ is injective, and 
 {\itshape unimodular} if $c$ is bijective. For
 $m\in\BQQ$, $L(m)$ denotes the same underlying group equipped with 
 the form multiplied by $m$, assuming that it is $\BZZ$-valued.
 $U,E_8$ and $D_4$ denote the lattices given by the 
 matrix $\left(
 \begin{smallmatrix}
 0 & 1 \\
 1 & 0
 \end{smallmatrix}
 \right)$,
 the Dynkin diagrams of type $E_8$ and $D_4$ respectively.
 We understand the latter two to be negative definite.\vspace{5mm}\\
 \quad A {\itshape finite quadratic form} is a triple $(A,q,b)$ 
 where $A$ is a finite abelian group, $q$ is a map $A\rightarrow 
 \BQQ/2\BZZ$ and $b$ is a bilinear map $A\times A\rightarrow \BQQ/\BZZ$ 
 which is symmetric and satisfies 
 \[q(x+y)= q(x)+q(y)+2b(x,y),\ x,y\in A. \] 
 \quad In the following we abbreviate $b(x,y)$
 (resp. $q(x)$) to $xy$ (resp. 
 $x^2$) and sometimes $(A,q,b)$ to $(A,q)$. We call $x^2$ the norm of $x$.
 As in the lattice case, we have a natural homomorphism 
 $c:A\rightarrow A^*=\Hom(A,\BQQ/\BZZ)$ defined by using $b$.
 $(A,q)$ is said to be {\itshape nondegenerate} if $c$ is bijective.
\vspace{5mm} \\
 For an even nondegenerate lattice $L$, we can canonically 
 associate a finite quadratic form $(A_L,q_L)$, called the discriminant
 quadratic form of $L$, by putting $A_L:= L^*/L$ and
 $q_L$ is the one naturally induced from the linear extension of the form
 on $L$ to $L^*\subset L\otimes \BQQ$.
 The discriminant group of $U(2)$ (resp. $D_4$) is denoted by 
 $u(2)$ (resp. $v(2)$).\\
 \quad For a lattice $L$, $O(L)$ (resp. $O(q_L)$) denotes the integral
 orthogonal group of $L$ (resp. of $(A_L,q_L)$). We note that 
 there is a natural homomorphism $\sigma_L:O(L)\rightarrow O(q_L)$.
 $L_{\BRR}$ (resp. $L_{\BCC}$) is the scalar extension of $L$ 
 to $\BRR$ (resp. $\BCC$).\vspace{5mm}\\

 The author is grateful to Professor Shigeru Mukai for many
 helpful discussions and suggestions. He indicated
 the example in Section \ref{example}. He is also grateful to Professors
 Shigefumi Mori and Noboru Nakayama for many valuable comments 
 throughout the seminars.
\section{Finiteness of conjugacy classes of finite
 subgroups}\label{finiteness}
 \quad First we collect some basic definitions about $K3$ surfaces.
 Let $X$ be a $K3$ surface. It is known that all $K3$ surfaces are 
 diffeomorphic. {\itshape A $K3$ lattice} is a lattice isomorphic to
 $H^2(X,\BZZ)=U^{\oplus 3}\oplus E_8^{\oplus 2}$. 
 $\omega_X$ is the period of $X$, namely $\BCC\omega_X=H^{2,0}(X)$.
 $NS(X)=\omega_X^{\perp}\cap H^2(X,\BZZ)$ is the N\'eron-Severi lattice of $X$.
 $T_X=(NS(X))^{\perp}$ is the transcendental lattice of $X.$\vspace{5mm}\\
 \quad We recall the structure of the integral automorphism group $O(NS)$ of $NS(X)$.
\begin{dfn}\quad
 \begin{enumerate}
  \item The {\it positive cone} $\mathcal{C}_X$ is the connected component of
	 $\{x\in NS(X)\otimes \BRR|\
        x^2>0\ \}$ which contains an ample divisor.
  \item The {\it ample cone} $\mathcal{A}_X$ is the subcone 
         of $\mathcal{C}_X$ generated as a semigroup by ample divisors 
	 multiplied by positive real numbers.
 \end{enumerate}
\end{dfn}
\begin{dfn}\quad 
 \begin{enumerate}
  \item The {\it Weyl group} $W_X$ of $X$ is the subgroup of $O(NS)$ generated by 
        automorphisms of the form $s_l: x\mapsto x+(xl)l$ for all elements 
	$l\in NS(X)$ with $l^2=-2$.
  \item $O^{\uparrow}(NS):=\{\vp\in O(NS)|\vp (\mathcal{C}_X)=\mathcal{C}_X\}$.
  \item $O^+(NS):=\{\vp \in O(NS)| \vp(\mathcal{A}_X)=\mathcal{A}_X\}$.
  \item $O_0(NS):=\ker (\sigma_{NS}: O(NS)\rightarrow O(q_{NS}))$.
 \end{enumerate}
\end{dfn}
 Further we use the abbreviations 
 \[ O_0^{\uparrow}(NS)=O_0(NS)\cap O^{\uparrow}(NS)\quad \text{\rm{and }} 
 \ O_0^{+}(NS)=O_0(NS)\cap O^{+}(NS).\]
 \quad $O^{\uparrow}(NS)$ is of index $2$ in $O(NS)$. The inclusions $O_0(NS)\subset
 O(NS)$ and $O^{+}_0(NS)\subset O^+ (NS)$ are of finite index since 
 $O(q_{NS})$ is a finite group.\\
 \quad The following relation between these subgroups are important.
\begin{prop}\label{dec}
 We have $(1)\ W_X\subset O_0^{\uparrow}(NS)$,
 $(2)\ O^{\uparrow}(NS)=W_X\rtimes O^{+}(NS)$ and 
 $(3)\ O_0^{\uparrow}(NS)=W_X\rtimes O_0^{+}(NS)$. 
\end{prop}
$Proof.$ Since the generator $s_l\in W_X$ acts trivially on the hyperplane $H_l$
 orthogonal to $l$ which intersects with the positive cone, $W_X$ preserves the 
 positive cone. If $x\in NS^{*}$, then $xl\in \BZZ$ and $s_l(x)=x$ modulo $\BZZ l$.
 This proves $(1)$. 
 We have the semidirect product decomposition as in $(2)$ because $W_X$ and $O^{\uparrow}(NS)$
 are discrete subgroups of 
 the isometry group of the Loba\u{c}evski\u{i}\ space
 modeled in $\mathcal{C}_X$ and $W_X$ is a 
 reflection group with ample cone as its fundamental domain. See \cite{Vinberg}.
 The proof of $(3)$ is the same. \qquad q.e.d.\\

 We proceed to the proof of Theorem \ref{FP}.
 For brevity, we say that a group $G$ has property $(FP)$ if 
 $G$ has only finitely many conjugacy classes of finite subgroups.
 For example, let $G$ be an algebraic group defined over $\BQQ$. Then 
 $G_{\BZZ}$ has property $(FP)$ by \cite[Theorem 4.3]{PR}
 which we call the theorem of Borel.
\begin{lem}\label{onFP}\quad \\
{\rm (1)} Let $\alpha : G \rightarrow K$ be a homomorphism of groups.
 If $\im \alpha$ has property (FP) and $\ker \alpha$ is finite, then $G$ has property (FP).\\
{\rm (2)} Let $G=W\rtimes K$ be a semidirect decomposition of a group $G$.
 If two finite subgroups $F_1,F_2\subset K$ are $G$-conjugate, then they are 
 also $K$-conjugate. In particular if $G$ has property (FP), then so does $K$.\\
{\rm (3)} Let $H\subset G$ be a subgroup of finite index. 
 If $G$ has property (FP), then so does $H$. 
\end{lem}
$Proof.$ (1) Let $P_1,\cdots,P_n$ be the complete representatives 
 of conjugacy classes of finite subgroups of $\im \alpha$. Then any conjugacy class of
 finite subgroups of $G$ has a representative included in at least one of $\alpha^{-1} (P_j),\ 
 (j=1,\cdots,n)$.\\
 (2) Assume $F_2=gF_1 g^{-1}$. $g$ can be written as $g=wk$, $w\in W,\ k\in K$.
 If $f_i\in F_i$ satisfy $f_2 =gf_1 g^{-1}$, then we 
 have $(k f_1k^{-1})^{-1}f_2 =(k f_1^{-1}k^{-1} wk f_1 k^{-1})w^{-1} \in W\cap K=\{ 1 \}$.
  Hence $kf_1k^{-1}=f_2$. Thus $F_1$ and $F_2$ are conjugate by $k\in K$.\\
 (3) Again let $P_1,\cdots,P_n$ be the complete representatives 
 of conjugacy classes of finite subgroups of $G$. We put $G/H=\{a_1H,\cdots,a_rH\}$.
 Then the conjugacy classes of finite subgroups of $H$ are represented by 
 $\{a_i^{-1}P_ja_i | i=1,\cdots,r\ {\rm and }\ j=1,\cdots, n\}$. 
 \qquad q.e.d.\vspace{5mm}\\

 Now we show Theorem \ref{FP}. In our words, 
\begin{thm}\label{finiteness}
Let $X$ be a $K3$ surface. Then $\Aut (X)$ has property $(FP)$.
\end{thm}
$Proof.$ First we assume that $X$ is projective.
 Consider the representation $r: \Aut(X)\rightarrow O(NS)$. Since every automorphism in 
 $\ker r$ fixes an ample divisor, $\ker r$ is finite. Thus it is enough 
 to show that $\im r$ has the property $(FP)$ by Lemma \ref{onFP}(1).\\
 \quad By the theorem of Borel above,
 $O(NS)$ has property $(FP)$. Then, by Lemma \ref{onFP}
 and Proposition \ref{dec}, the property $(FP)$ goes down to $O^{\uparrow}(NS)$ and 
 $O^+(NS)$. Now by the global Torelli theorem
 \cite{PSS}, $\im r$ contains $O^+_0(NS)$, since 
 $\vp\in O^+_0(NS)$ preserves the ample cone $\mathcal{A}_X$ and can be extended 
 to an isometry of $H^2(X,\BZZ)$ which acts trivially on $T_X$.
 Thus we obtain \[O^+_0(NS)\subset \im r \subset O^+(NS)\]
 and these inclusions are of finite index. 
 We see that $\im r $ has the property $(FP)$ by Lemma \ref{onFP}(3).\\
 \quad Next we assume that $X$ is not projective.
 Nikulin \cite{Nikulin2} shows that any automorphism of $X$ of finite order 
 acts on $T_X$ trivially. Therefore it is enough to show that
 $G=\ker (\Aut(X) \rightarrow O(T_X))$ has property $(FP)$.
 We consider the representation $r: G\rightarrow O(NS)$.\\
 \quad If $\algdim (X)=0$, then $NS(X)$ is negative definite. 
 Hence $NS(X)\oplus T_X$ is of finite index in $H^2(X,\BZZ)$, $r$ is injective 
 and $G$ is finite since $O(NS)$ is a finite group.\\
 \quad If $\algdim (X)=1$, then $NS(X)$ has one-dimensional kernel $\BZZ e$ and 
 $Q:=NS(X)/\BZZ e$ is negative definite. Every element of $G$ fixes $e$ since 
 $e^2=0$ and exactly one of $e$ and $-e$ is represented by an effective cycle.
 Thus $r$ induces $s: G\rightarrow O(Q)$. Let $g\in G$ be an element of finite order.
 Since the fixed part $H^2(X,\BZZ)^g$ is nondegenerate by the lemma below,
 it follows that if $s(g)=\id_Q$ then $g=\id_X$.
 On the other hand, $O(Q)$ is a finite group. Thus $\Aut(X)$ has only
 finitely many elements of finite order and $\Aut(X)$ has property $(FP)$.  
 \qquad q.e.d.\
 \begin{lem}\label{fff}
 \quad Let $L$ be a nondegenerate lattice and let $g$ be an isometry of $L$ of 
 finite order $n$. Let $M=L^g =\{x\in L| gx=x\}$ be the fixed lattice.
 Then $M$ is nondegenerate.
\end{lem}
 $Proof.$ Let $0\neq x\in M$. Since $L$ is nondegenerate, there 
 exists $y\in L$ with $xy\neq 0$. Put $z=y+g(y)+\cdots+g^{n-1}(y)$.
 Obviously $z\in M$ and we have 
 \[ xz=xy+g(x)g(y)+\cdots+g^{n-1}(x)g^{n-1}(y)=nxy\neq 0.\]
 Therefore $M$ is nondegenerate. \qquad q.e.d. \vspace{5mm}\\
 \quad As a corollary, Corollary \ref{fEq} follows.
 In fact, if two free involutions $i_1$ and $i_2$ are conjugate by 
 an automorphism $g$, then $g$ induces an isomorphism between
 $X/i_1$ and $X/i_2$.
\section{Number of distinct Enriques quotients}\label{number}
 \quad The isomorphism classes of Enriques quotients are exactly the 
 conjugacy classes of free involutions by the next proposition.
\begin{prop}\label{lifting}
 \quad Let $X$ be a $K3$ surface and let $i_1$ and $i_2$ be
 free involutions on $X$. Then, $X/i_1$ and $X/i_2$ are
 isomorphic if and only if there exists
  an automorphism $g$ of $X$ such that $g i_1 g^{-1}=i_2$.
\end{prop}
 $Proof.$ The ``if'' part is trivial; see the sentence after 
 Lemma \ref{fff}. Conversely, let $h$ be an isomorphism
 from $Y_1:=X/i_1$ to $Y_2:=X/i_2$.
 It induces an isomorphism of the canonical line bundles
 $h^* \mathcal{K}_{Y_2} \rightarrow \mathcal{K}_{Y_1}$.
 Since $X$ and $Y_1,Y_2$ are related as $X=
 \Spec(\mathcal{O}_{Y_j}\oplus
 \mathcal{K}_{Y_j})\  j=1,2$, $h$ induces the 
 following commutative diagram\\
\[
\begin{CD}
  X          @>g>>         X       \\
  @V{\pi}VV            @VV{\pi}V   \\
  X/{i_1}    @>>h>         X/{i_2}
\end{CD}
\]
 where $\pi$ denotes the covering map. It is clear that $g$ is 
 the desired automorphism of $X$. \qquad q.e.d.\\
 \quad \\
 \quad We put 
 \[\mathfrak{M}:= \left\{ 
 \begin{array}{c|c}
   M\subset NS & \begin{array}{l}
               \text{\rm{a primitive sublattice which satisfies }}\\
	       \text{\rm (A)}:\ M\cong U(2)\oplus E_8(2)\\
	       \text{\rm (B)}:
		\text{\rm{ No vector of square $-2$ in $NS(X)$ is orthogonal to $M$.}}
	       \end{array}
 \end{array}
 \right\}. \]
 \quad Recall that $U\oplus E_8$ is the Enriques lattice (modulo torsion) and 
 $U(2)\oplus E_8(2)$ is the pullback in the covering $K3$ lattice. 
 For each $M\in \mathfrak{M}$, we define an 
 isometry $i_M : H^2(X,\BZZ) \rightarrow H^2(X,\BZZ)$ by 
 $i_M(m)=m$ when $m\in M$ and $i_M(n)=-n$ when $n$ is orthogonal to $M$.
 This is well-defined because $M\cong U(2)\oplus E_8(2)$ is $2$-elementary.\\
\begin{prop}\label{correspondence} 
 \quad On a $K3$ surface $X$, there is a one-to-one correspondence 
 between free involutions on $X$ and primitive sublattices $M$
 of $NS(X)$ which satisfy {\rm (A)} and {\rm (B)} above and the following\\
  \qquad {\rm (C)}\ $M$ contains an ample divisor. \\
 \quad In other words, $i_M$ defined above is a free involution 
 if and only if $M$ contains an ample divisor. Also any free involution can 
 be written in the form $i_M$.
\end{prop}
 $Proof.$ We associate a free involution with its invariant sublattice in 
 $H^2(X,\BZZ)$. The statement follows from \cite[Corollary 2.5]
 {Namikawa}, \cite[Theorem 4.2.2, p1426]{Nikulin3} and the strong 
 Torelli theorem for $K3$ surfaces \cite{PSS}. 
 In \cite{Nikulin3}, the assumption is slightly  
 different from ours, but the same proof goes.\qquad q.e.d.\\
 \quad \\
 \quad To count the number of distinct Enriques quotients, we consider the 
 natural action of $O(NS)$ on $\mathfrak{M}$,
 \[ O(NS)\ni\vp : M\mapsto \vp(M)\in \mathfrak{M}.\]
 Corresponding lattice automorphisms satisfy $i_{\vp(M)}=\vp i_M\vp^{-1}$.\\
 \quad In the following, $\Aut (T_X,\omega_X)$ is the 
 subgroup of $O(T_X)$ consisting of the 
 integral orthogonal transformations
 whose scalar extention to $\BCC$ preserves 
 the period $\BCC\omega_X\subset T_X\otimes \BCC.$ 
\begin{thm}\label{bound}
 \quad Let $M_1,\cdots,M_k\in \mathfrak{M}$ be a (finite) complete
 set of representatives for the action of $O(NS)$ on $\mathfrak{M}$.
 For each $j=1,\cdots,k$, let
 \[ K^{(j)}=\{\vp\in O(NS)| \vp(M_j)=M_j \}\]
 be the stabilizer subgroup of $M_j$ and $\sigma(K^{(j)})$ its
 canonical image in $O(q_{NS})$. We put
 \[B_0=\sum^k_{j=1} \#(O(q_{NS})/\sigma(K^{(j)})).\]
 \quad {\rm (1)}\ The number of distinct Enriques quotients of $X$
 does not exceed $B_0$.\\
 \quad {\rm (2)}\ If $\sigma:O(NS)\rightarrow
 O(q_{NS}) $ is surjective and if $\Aut(T_X,\omega_X)=\{\pm\id\}$,
 then $X$ has exactly $B_0$ distinct Enriques quotients.  
 \end{thm}
$Proof.$
 First we remark that by Proposition 1.15.1 in \cite{Nikulin1},
 the set of representatives is always finite.
 In view of Proposition \ref{lifting}, we can count the number 
 of distinct Enriques quotients separately for each orbit $O(NS)\cdot M_j$.
 Hence, for simplicity, we fix an orbit and omit the index $j$ 
 so that we use the symbols $M:=M_j,\ \mathcal{O}:=O(NS)\cdot M$ 
 and $K=K^{(j)}$.\\
\stepone $\mathcal{O}$ contains an element which corresponds to a 
 free involution.\\
 \quad $Proof.$ The following is a standard argument used in \cite{PSS}. 
 Our proof is taken from \cite{Namikawa}. 
 By the condition (A) in Proposition \ref{correspondence},
 $M\cap \mathcal{C}_X\neq \emptyset$.
 Consider in $\mathcal{C}_X$ 
 countably many hyperplanes $H_d=\{x\in NS_{\BRR}| xd=0\}$, where
 $d$ runs over $(-2)$ vectors in $NS$. The union $\cup H_d$ is a 
 locally finite closed subset in $\mathcal{C}_X$ and does not contain
 $M$ by the condition (B).
 The complement $\mathcal{C}_X -\cup H_d$ is a collection of 
 (at most) countably many connected 
 open sets, namely chambers, which corresponds to the elements of $W_X$
 in one-to-one way. The ample cone $\mathcal{A}_X$ 
 equals one of the chambers. Thus if we choose 
 $v\in M \cap\mathcal{C}_X -\cup H_d$,
 there exists $\vp \in W_X \subset O(NS)$ 
 such that $\vp (v)$ is an ample divisor. $\square$\\
 \quad Thus we can assume that $i_M$ is already a free involution of $X$. Next we set
 \[\mathfrak{N}:=\{M'\in \mathcal{O}|
 i_{M'}\text{\rm{ is a free involution}}\}.\]
\steptwo $\mathfrak{N}=O^+(NS)\cdot M$.\\
 \quad $Proof.$ $\supset$ follows from Proposition \ref{correspondence}.
 Let $\vp \in O(NS)$ and 
 suppose $i_{\vp(M)}$ is a free involution.
 We can assume $\vp\in O^{\uparrow}(NS)$, since otherwise 
 $-\vp\in O^{\uparrow}(NS)$ and $(-\vp)(M)=\vp(M)$.
 By Proposition \ref{dec}, we can write $\vp=w \psi$
 with $w\in W_X$ and $\psi \in O^+(NS)$. Using Lemma 
 \ref{onFP}(2), $i_{\vp(M)}=w i_{\psi(M)} w^{-1}$ 
 implies $i_{\vp(M)}=i_{\psi(M)}$. Therefore $\vp(M)=\psi(M)$. $\square$\\
 \stepthree Let $\psi_1, \psi_2 \in O^{+}(NS)$. If $\sigma(\psi_j )$
 have the same class in $O(q_{NS})/\sigma(K)$, then $i_{\psi_1 (M)}$ and 
 $i_{\psi_2 (M)}$ are conjugate in $\Aut (X)$.\\
 \quad $Proof.$ By the assumption
 $\sigma(\psi_1^{-1}\psi_2)\in \sigma(K)$, so 
 there exists $\vp\in K$ such that $\sigma(\vp)=\sigma(\psi_1^{-1}\psi_2)$.
 It follows that $\sigma(\psi_1\vp\psi_2^{-1})=\id$, 
 so that $\psi_1\vp\psi_2^{-1}|_{NS}$ together with $\id_{T_X}$ gives 
 an automorphism $a$ of $X$, by the Torelli theorem. It follows that 
 \[ ai_{\psi_2(M)}a^{-1} = i_{a\psi_2(M)}= i_{\psi_1 \vp(M)}
 =i_{\psi_1(M)}.\square\]
 \quad By now, we have proved that $\mathcal{O}$ contains at most
 $\# O(q_{NS})/\sigma(K)$ distinct Enriques quotients.
 Assertion (1) follows.\\
 \stepfour If $\Aut(T_X,\omega_X)=\{\pm\id\}$, then the converse
 of Step 3 holds.\\
 \quad $Proof.$ 
 Assume there exists $\vp\in \Aut(X)$ such that 
 $\vp i_{\psi_1(M)} \vp^{-1}= i_{\psi_2(M)}$, which is equivalent to 
 $\vp\psi_1(M)=\psi_2(M)$ and to $\psi_2^{-1} \vp \psi_1\in K$.
 By the assumption, $\sigma(\vp)= \pm \id$ which is contained in the 
 center of $O(q_{NS})$. It follows that 
 $\pm \sigma (\psi_2^{-1} \psi_1)\in \sigma(K).$
 We remark that $-\id\in \sigma(K)$ since $\sigma(i_M)=-\id$.
 Therefore we get $\sigma(\psi_2^{-1}\psi_1)\in \sigma(K)$.$\square$\\ 
 \stepfive If $\sigma$ is surjective, then the restriction
  $\sigma|_{O^+(NS)}$ is also surjective.\\
 \quad $Proof.$
 Put $N=M^{\perp}$ in $NS$.
 Since $M$ is $2$-elementary, $t=(-\id_M,\id_N)$ extends to an isometry 
 of $NS$. $t$ doesn't preserve the positive cone. Therefore
 $O(NS)$ is generated by 
 $t$ and $O^{\uparrow}(NS)$. This implies the surjectivity of 
 $O^{\uparrow}(NS)\rightarrow O(q_{NS})$. 
 By Proposition \ref{dec}, the assertion follows.$\square$\\
 \quad Now the proof is complete. \qquad q.e.d.\\
 \quad \\
 \quad Lastly we mention a useful theorem of Nikulin
 in \cite{Nikulin1}
 which saves us from checking one of the conditions in 
 Theorem \ref{bound} (2).
\begin{thm}[Nikulin]\label{useful}
 Let $T$ be an even indefinite nondegenerate lattice satisfying the 
 following two conditions:\\
 \quad {\rm (1)}\ $\rank(T)\ge l(A_{T_p})+2$ for all prime numbers $p$ 
 except for $2$.\\
 \quad {\rm (2)}\ if $\rank(T)=l(A_{T_2})$, then $q_{T_2}$ contains a  
 component $u(2)$ or $v(2)$.\\
 \quad Then the genus of $T$ contains only one class, 
 and the homomorphism $O(T)\rightarrow O(q_T)$ is surjective.
 Here $A_{T_p}$ denotes the $p$-component of the finite abelian group
 $A_T$ and $l$ denotes the number of minimal generators. $\square$
\end{thm}
\section{Enriques quotients of $K3$ surfaces in the Heegner divisors}\label{unboundedness}
 \quad In this section we prove Theorem \ref{unbound}.
 We deal with certain divisors of the period domain $\MFD$
 of $U(2)\oplus E_8(2)$-polarized marked 
 $K3$ surfaces.
 Fix the unique primitive embedding of $U(2)\oplus E_8(2)$ in the $K3$ lattice 
 $\Lambda$. Then $\MFD$ is by definition 
\[ \MFD
 := \{[\omega]\in \BPP((U(2)\oplus E_8(2))^{\perp}_{\BCC})| 
 \omega^2 =0,\ \omega\overline{\omega}>0\}. \]
 Here $\BPP(V)$ means the associated projective space of
 a complex vector space $V$, which consists of all lines 
 through the origin. It follows from the surjectivity of the period map 
 that every point of $\MFD$ corresponds to a $K3$ surface $X$ with a marking 
 $H^2(X,\BZZ)\cong \Lambda$.\\
 \quad Let $S\subset \Lambda$ be a primitive sublattice of rank $11$ containing 
 the lattice $U(2)\oplus E_8(2)$ fixed above.
 Then the subset 
\[ \MFD (S)
 := \{[\omega]\in \BPP(S^{\perp}_{\BCC})| 
 \omega^2 =0,\ \omega\overline{\omega}>0\} \]
 is called {\it the Heegner divisor of type} $S$ in $\MFD$.
 Let $X$ be a marked $K3$ surface whose period is in $\MFD (S)$.
 Since $NS(X)$ is written as $\Lambda\cap \omega_X^{\perp}$,
 $NS(X)$ contains the primitive
 sublattice $S$.
\begin{prop}\label{generic}
 If $X$ corresponds to a very general point of $\MFD (S)$, namely 
 to a point in the complement of a union of countably many closed 
 analytic subset of $\MFD (S)$, then we have $NS(X)=S$
 and $\Aut (T_X,\omega_X)=\{\pm \id \}$.
\end{prop} 
$Proof.$ This is a well-known fact. For the latter, the same 
 proof as in \cite[Lemma 2.9]{B-P} works.\qquad q.e.d.\\
 \quad \\
 \quad We consider the case when 
 \[S=U(2)\oplus E_8(2)\oplus \langle -2N \rangle,\]
 where $\langle -2N \rangle$ is the rank $1$
 lattice whose generator $g$ has $g^2=-2N$.
 It is easy to see that the $K3$ lattice $\Lambda$ contains $S$ as a 
 primitive sublattice. We fix it once and for all.
 The discriminant form of $S$ is 
 isomorphic to $q=u(2)^{\oplus 5}\oplus c(-2N)$, where $c(-2N)$ is the discriminant 
 form of $\langle -2N \rangle$.
 
 Let the integer $N$ be $N=4p_1 \cdots p_l$, where $p_1,\cdots, p_l$
 are distinct odd prime numbers.
 In the next we compute the order of $O(q)$.
\begin{lem}
 \[\#\{ x\in c(-2N)|\ord (x)=2N, x^2\equiv -1/2N\
 (\mymod 2\BZZ)\} =2^{l+1},\]
 \[\#\{ x\in c(-2N)|\ord (x)=2N, x^2\equiv 1-1/2N\ (\mymod 2\BZZ)\}
 =2^{l+1}.\]
\end{lem}
$Proof.$ The left-hand-side of the first equality is 
\[
\begin{array}{ll}
\  &  \#\{ k\in \BZZ| (k,2N)=1,1\le k\le 2N-1\ 
 \text{\rm  and } -k^2/2N=-1/2N \in \mathbb{Q}/2\BZZ\}\\
= &  (1/2)\#\{ k\in \BZZ| (k,4N)=1,1\le k\le 4N-1 \ \text{\rm  and }
 k^2-1\equiv 0 (\mymod 4N) \}\\
= &  (1/2)\#\{ x\in (\BZZ/4N\BZZ)^{\times}| \text{ord}(x)=1\ 
 \text{\rm  or }2\}.
\end{array}
\]
 \quad Then we can use the structure theorem of the unit group 
$(\BZZ/4N\BZZ)^{\times}$. Similarly the left-hand-side of the latter is 
\[
\begin{array}{ll}
\  &  \#\{ k\in \BZZ| (k,2N)=1,1\le k\le 2N-1
 \ \text{\rm  and }-k^2/2N=1-1/2N \in \mathbb{Q}/2\BZZ\}\\
= &  (1/2)\#\{ k\in \BZZ| (k,4N)=1,1\le k\le 4N-1
 \ \text{\rm  and }k^2\equiv 1-2N (\mymod 4N)\}\\
= &  (1/2)\#\{ k\in \BZZ| (k,4N)=1,1\le k\le 4N-1,k^2\equiv 1 (\mymod 2N)
    \ \text{\rm  and }k^2\not\equiv 1 (\mymod 4N)\}.\\
\end{array}
\]
\quad Using the commutative diagram
\[
\begin{CD}
(\BZZ/4N\BZZ)^{\times} @>{\sim}>> 
  (\BZZ/2^{4}\BZZ)^{\times}  @. \oplus 
  (\BZZ/p_1\cdots p_l\BZZ)^{\times} \\
@V{\alpha}VV          @V{\beta}VV    @V{\sim}VV  \\
(\BZZ/2N\BZZ)^{\times} @>{\sim}>> 
  (\BZZ/2^3\BZZ)^{\times}  @. \oplus 
  (\BZZ/p_1\cdots p_l\BZZ)^{\times} \\
\end{CD}
\]
where $\alpha,\ \beta$ are both 2:1 maps, we can count the number of  
elements which have order 2 
in the bottom row but do not in the top row. \qquad q.e.d. \\
\begin{prop}\label{auto}
 $O(q)$ acts transitively on the set of elements $x\in q$ 
 with $x^2\equiv -1/2N\ (\mymod 2\BZZ)$. There are $2^{11+l}$ such elements.
\end{prop}
$Proof.$
 Such element $x$ generates a subgroup $\langle x\rangle$ isomorphic to $c(-2N)$.
 Since it is nondegenerate, $\langle x \rangle$ is a direct summand in $q$.
 This implies the transitivity.
 If we put the number of elements in $u(2)^{\oplus 5}$ with norm $0$ to be $A$
 and norm $1$ $2^{10}-A$, we can compute the 
 length of the orbit 
 as $2^{l+1}\cdot A+2^{l+1}\cdot (2^{10}-A)=2^{11+l}$. \qquad q.e.d.\\
 \quad \\
 \quad Under these computations we can prove
\begin{thm}
 For any nonnegative integer $l$, there exists a $K3$ surface $X$ with 
 exactly $2^{l+10}$ distinct Enriques quotients.
\end{thm}
$Proof.$ By Proposition \ref{generic}, there exists a $K3$ surface $X$ 
 such that $NS(X)\cong S$ and $\Aut(T_X,\omega_X)=\{\pm \id\}$.
 We show that the primitive embedding of
 $U(2)\oplus E_8(2)$ into $NS(X)$ is unique under the action of $O(NS)$.
 In fact, since $NS(1/2)$ is again an even lattice, we have
 a natural identification 
\[ \Hom(U(2)\oplus E_8(2), NS)=\Hom(U\oplus E_8, NS(1/2)).\]
 We see that any primitive embedding as above is a direct 
 summand. This clearly implies the uniqueness.\\
 \quad Obviously $NS(X)$ has a primitive sublattice $M$ isomorphic 
 to $U(2)\oplus E_8(2)$ and $M^{\perp}=\langle -2N \rangle$. Let
 $K$ be the stabilizer group of $M$ and $\sigma(K)$ its canonical 
 image in $O(q_{NS})$. Since $NS\cong M\oplus M^{\perp}$ we see that 
 $K=O(M)\times O(M^{\perp})=O(M)\times \{\pm\id_{M^{\perp}}\}$.
 On the other hand by Theorem \ref{useful} $\sigma_M:O(M)\rightarrow
 O(q_M)$ is surjective. This shows $\sigma(K)=O(u(2)^{\oplus 5})
 \times \{\pm \id\}\subset O(q_M\oplus q_{M^{\perp}})=O(q_{NS})$.
 Thus $\#(O(q_{NS})/\sigma(K))=\#O(q_{NS})/2\#O(u(2)^{\oplus 5})
 =2^{10+l}$ by Proposition \ref{auto}. This together with
 Theorem \ref{bound},\ref{useful} completes the proof.
 \qquad q.e.d.  \quad \vspace{5mm}\\
 In fact we can classify all the possible N\'eron-Severi lattices
 of a $K3$ surface with Picard number $11$ having an Enriques quotient.
 In each case, we can compute the number of Enriques quotients as follows
 by an explicit calculation. Details are omitted.
 The result is as follows.
\begin{prop}
 Let $X$ be a $K3$ surface with Picard number 
 $11$ having an Enriques quotient. Then the N\'eron-Severi lattice of $X$
 is one of the followings.
\[\begin{array}{lc}
  \text{\rm Type I}: &
  U(2)\oplus E_8(2)\oplus \langle -2N\rangle \qquad (N\ge 2)\\
  \text{\rm Type II}: &
  U\oplus E_8(2)\oplus \langle -4M\rangle \qquad (M\ge 1).
\end{array} 
\]
\end{prop}
 If we put $2N=2^e p_1^{e_1} \cdots p_l^{e_l}$ in type I, or 
 $4M=2^e p_1^{e_1} \cdots p_l^{e_l}$ in type II, the bound $B_0$ 
 in Theorem \ref{bound} is as follows.
\begin{prop}\label{div}
 \quad $B_0=
\begin{cases}
  2^{l-1} & \text{in Type I and } e=1\\ 
  (2^5+1)\cdot2^{l+4} & \text{in Type I and } e=2 \\
  2^{l+10} & \text{in Type I and } e \ge 3 \\
  1 & \text{in Type II and } e=2, l=0 \\
  2^{l-1} & \text{in Type II and } e=2, l>0 \\
  2^{2l+5}& \text{in Type II and } e \ge 3 
\end{cases}
$\end{prop}
 The lattice $S$ we used fits in the third case.
\section{Enriques quotients of generic Kummer
 surfaces of product type}\label{example}
\begin{figure}[htbp]
 \begin{center}
 \psfrag{$G_{11}$}{$G_{11}$}
 \psfrag{$G_{12}$}{$G_{12}$}
 \psfrag{$G_{13}$}{$G_{13}$}
 \psfrag{$G_{14}$}{$G_{14}$}
 \psfrag{$G_{21}$}{$G_{21}$}
 \psfrag{$G_{31}$}{$G_{31}$}
 \psfrag{$G_{41}$}{$G_{41}$}
 \psfrag{E1}{$E_1$}
 \psfrag{E2}{$E_2$}
 \psfrag{E3}{$E_3$}
 \psfrag{E4}{$E_4$}
 \psfrag{F1}{$F_1$}
 \psfrag{F2}{$F_2$}
 \psfrag{F3}{$F_3$}
 \psfrag{F4}{$F_4$}
\includegraphics{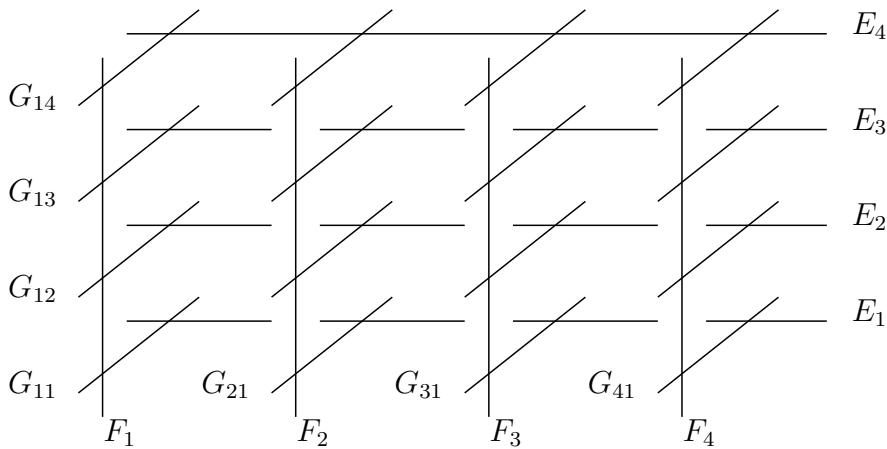}
 \end{center}
 \caption{the double Kummer configuration}
 \label{fig}
\end{figure}

\begin{figure}[htbp]
 \begin{center}
 \psfrag{e1}{$e_1$}
 \psfrag{e2}{$e_2$}
 \psfrag{e3}{$e_3$}
 \psfrag{e4}{$e_4$}
 \psfrag{e5}{$e_5$}
 \psfrag{e6}{$e_6$}
 \psfrag{e7}{$e_7$}
 \psfrag{e8}{$e_8$}
\includegraphics[height=2cm,clip]{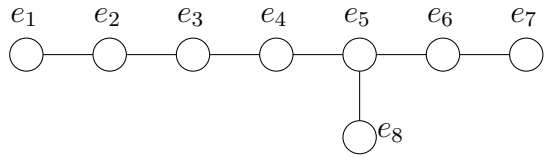}
 \end{center}
 \caption{Dynkin diagram of type $E_8$}
 \label{fig2}
\end{figure}

\begin{figure}[htbp]
 \begin{center}
 \psfrag{C1}{$C_1$}
 \psfrag{C2}{$C_2$}
 \psfrag{1}{$e_1\!\! +\!\! f_1\!\! +\!\! e_2\!\! +\!\! f_2$}
 \psfrag{2}{$e_1\!\!+\!\!f_2$}
 \psfrag{3}{$e_2\!\!+\!\!f_1$}
 \psfrag{4}{$e_1\!\!+\!\!e_2$}
 \psfrag{5}{$e_1$}
 \psfrag{6}{$e_2$}
 \psfrag{7}{$f_1\!\!+\!\!f_2$}
 \psfrag{8}{$f_2$}
 \psfrag{9}{$f_1$}
\includegraphics[height=5cm,clip]{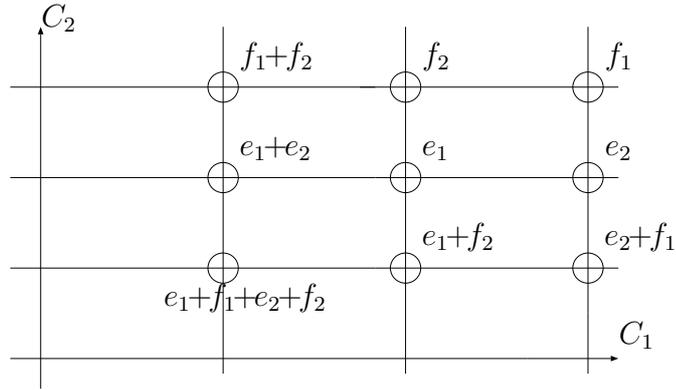}
 \end{center}
 \caption{Lieberman involutions}
 \label{fig3}
\end{figure}

\begin{figure}[htbp]
 \begin{center}
 \psfrag{1}{$e_2\! +\! f_2$}
 \psfrag{2}{$f_1\!+\!e_2\!+\!f_2$}
 \psfrag{3}{$e_1\!+\!f_1\!+\!e_2$}
 \psfrag{4}{$e_1\!+\!e_2\!+\!f_2$}
 \psfrag{5}{$e_1\!+\!f_1$}
 \psfrag{6}{$e_1\!+\!f_1\!+\!f_2$}
\includegraphics[height=11cm,clip]{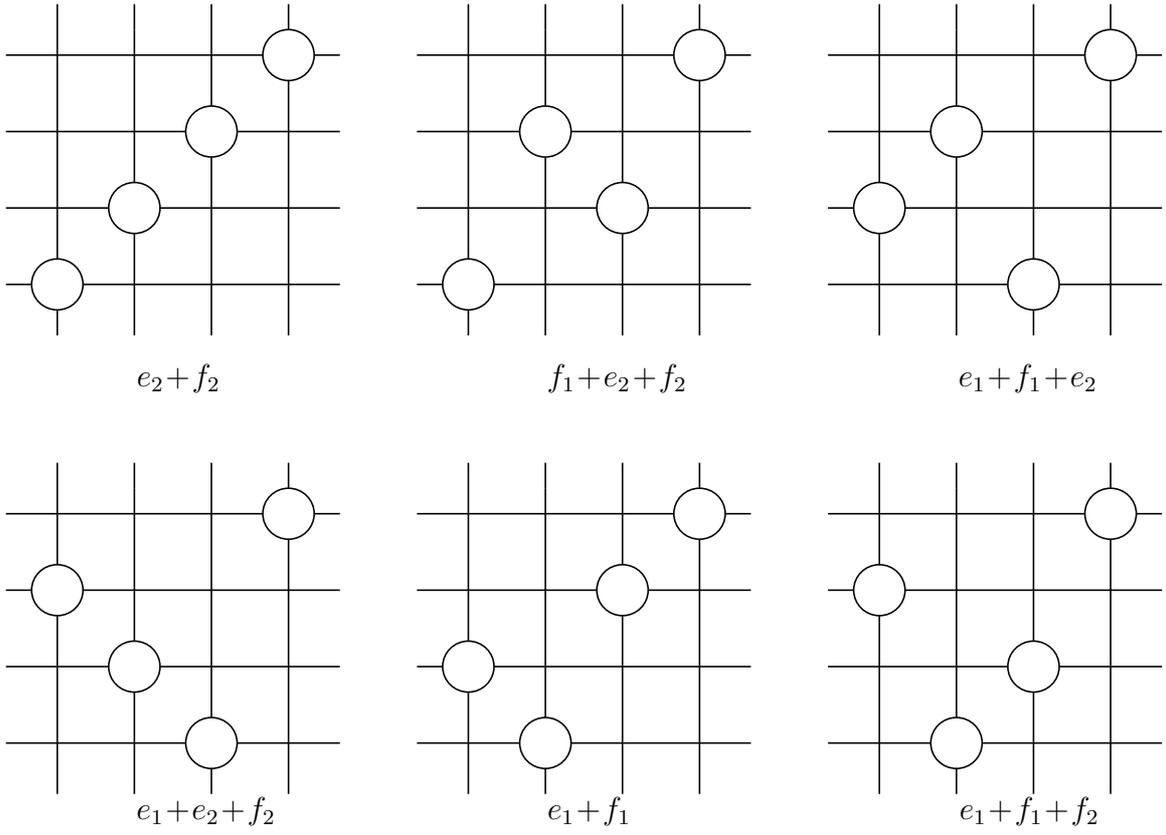}
 \end{center}
 \caption{Kondo-Mukai involutions}
 \label{fig4}
\end{figure}

 \quad Kondo found the first example 
 of a $K3$ surface which has two distinct Enriques quotients
 in \cite[Remark 3.5.3]{Kondo}, where 
 he computed the automorphism groups of the two quotients.
 Recently Mukai generalized Kondo's example which we now describe.\\
 \quad \\
 \quad {\bf Kummer surfaces of product type.}
 Let $C_1$ and $C_2$ be elliptic curves
 and construct the Kummer surface as $X=Km(C_1 \times C_2)$.
 We put the $2$-torsion points of $C_1$ (resp. $C_2$) as $\{b_1=0,
 b_2, b_3, b_4\}$ (resp. $\{c_1=0, c_2, c_3, c_4\}$) and denote 
 by $\delta$ the natural rational map of $C_1\times C_2$ to $X$.
 Let $E_k$ (resp. $F_k$) be the image of $C_1\times \{c_k\}$
 (resp. $\{b_k\}\times C_2$) by $\delta$. Then $X$ has the configuration
 of $24$ smooth rational curves as in Figure \ref{fig}, where $G_{ij}$ is 
 the exceptional curve corresponding to $(b_i,c_j)\in C_1 \times C_2$.\\
 \quad Sometimes it is called the double Kummer configuration. 
 In the following we introduce two kinds of free involutions on $X$
 with parameters.\vspace{3mm}\\
 \quad {\bf Involutions of Lieberman type.}
 Let $a=(b_i,c_j)\in C_1\times C_2$ be a $2$-torsion point not lying on 
 $C_1\times \{ 0 \}$ or $\{0\}\times C_2$. Let $\tau$ (resp. $\rho$) be 
 the involution of $X$ induced by the automorphism
 $(-\id_{C_1},\id_{C_2})$ (resp. the translation by $a$)
 of $C_1\times C_2$. Then $\sigma=\tau\rho$ is a free involution.
 We remark that $a$ has the parameter $i,j$ with $2\le i,j \le 4$.\vspace{2mm}\\
 \quad {\bf Involutions of Kondo-Mukai type.}
  Let $\tau$ be as above. $X/\tau$ is isomorphic to 
 $\BPP^1\times \BPP^1$ with sixteen points blown up, which correspond
 to $G_{ij}$'s in Figure \ref{fig}. We regard $\BPP^1\times \BPP^1$
 as a smooth quadric in $\BPP^3$ so that there are morphisms 
\[\begin{CD}
 X @>>> X/\tau @>>{\varepsilon}> \BPP^1\times \BPP^1\subset \BPP^3.
\end{CD}\]
 \quad Choose two permutations $I=\{i_1,i_2,i_3,i_4\}$ and
 $J=\{j_1,j_2,j_3,j_4\}$ of $\{1,2,3,4\}$ and we put 
 $g_k=\varepsilon(G_{i_k j_k})\in \BPP^1\times \BPP^1$. We project
 $\BPP^1\times \BPP^1$ onto a hyperplane $\cong \BPP^2$ of $\BPP^3$ from $g_4$.
 This birational map $\varepsilon'$ contracts two curves
 whose images we denote by $P$ and $Q \in \BPP^2$.
 Let $\overline{g}_k = \varepsilon' (g_k),\ k=1,2,3$.
 Then we obtain two involutions on $X$:
 One is the covering transformation $\nu$
 of the degree $2$ map $X\rightarrow \BPP^2$.
 The other involution $\mu$ is induced by the unique Cremona transformation of 
 $\BPP^2$ centered at $\overline{g}_1,\overline{g}_2,\overline{g}_3$
 interchanging $P$ and $Q$.
 Then $\sigma=\nu\mu$ is a free involution.
 It can be shown that $\sigma$ 
 depends only on the set $\{G_{i_k j_k}\}_{1\leq k\leq 4}$.
 See \cite{Mukai} for such canonical description of this involution.
 Therefore the parameter is the set $\{G_{i_k j_k}\}_{1\leq k\leq 4}$.
 There are twenty-four choices of 
 parameters in appearance, but 
 we will reveal that there are essentially six. $\square$\\ 
 \quad \\
 \quad If $C_1$ and $C_2$ are chosen to be general enough, we see that 
 \[ NS(X)\cong U\oplus E_8\oplus D_4^{\oplus 2},\ 
 T_X\cong U(2)^{\oplus 2},\ A_{NS(X)}=u(2)^{\oplus 2}\quad 
 {\rm and}\quad  \Aut(T_X,\omega_X)=\{\pm\id\}.\]
 \quad We will call such $X$ {\it a generic Kummer surface of product type.}
 In the following main theorem of this section 
 we classify the all free involutions on $X$.
\begin{thm}\label{Kummer}
 \quad Let $X$ be a generic Kummer surface of product type.
 Then $X$ has exactly fifteen distinct Enriques quotients which 
 are naturally in one-to-one correspondence with nonzero elements
 of $A_{NS(X)}$. Moreover all of them can be geometrically 
 constructed from the preceding examples by choosing appropriate parameters.
\end{thm}
 \quad We remark that the Lieberman involutions correspond to nine 
 elements of norm $0$ of $A_{NS(X)}$ and Kondo-Mukai involutions 
 to six elements of norm $1$.
 \quad In the rest of this section we prove 
 Theorem \ref{Kummer}. First we determine the isomorphism classes
 of primitive embeddings of $M:=U(2)\oplus E_8(2)\subset
 S=U\oplus E_8\oplus D_4^{\oplus 2}.$ This step is purely 
 lattice-theoretic.
 We use the following theorem of Nikulin \cite[Propositions 1.5.1 and 1.15.1]{Nikulin1}.
\begin{thm}\label{lemb}
 Let $M\subset S$ be a primitive embedding of even nondegenerate lattices.
 Put $N:=M^{\perp}$.
 Then the following isomorphisms exist.
 (Note that {\rm (1)} expresses $A_S$ in terms of $A_M$ and $A_N$ while {\rm (2)} expresses 
 $A_N$ in terms of $A_M$ and $A_S$.)\\
 {\rm (1)} There are subgroups $\Gamma_M\subset A_M,\Gamma_N\subset A_N$ 
 and a sign-reversing isometry $\gamma: \Gamma_M\rightarrow \Gamma_N$
 such that if $\Gamma$ is the pushout of $\gamma$, namely 
 $\Gamma=\{ (x,\gamma(x))\in A_M\oplus A_N | x\in \Gamma_M\}$,
 then \[ q_S\cong (q_M\oplus q_N|_{\Gamma^{\perp}}/\Gamma).\]
 {\rm (2)} There are subgroups $\Gamma_M\subset A_M,\Gamma_S\subset A_S$ 
 and a sign-reversing isometry $\gamma: \Gamma_M\rightarrow \Gamma_S$
 such that if $\Gamma$ is the pushout of $\gamma$, 
 then \[q_N\cong (-q_M\oplus q_S|_{\Gamma^{\perp}}/\Gamma).\]
\end{thm} 
\begin{prop}\label{E82}
 Let
 \[S=U\oplus E_8\oplus D_4^{\oplus 2}\ {\rm and}\ M=U(2)\oplus E_8(2).\]
 Assume $M$ is a primitive sublattice of $S$ which is orthogonal to 
 no $(-2)$ vectors of $S$. Then\\
 {\rm (1)} $N:=M^{\perp}$ is isomorphic to $E_8(2)$.\\
 {\rm (2)} There are exactly two such primitive sublattices up to 
 the action of $O(S)$.
\end{prop}
$Proof.$ By Theorem \ref{lemb} (2), $A_N$ is a $2$-elementary abelian group
 and $q_N$ takes only integral values on $A_N$.
 On the other hand, $N$ is a negative definite lattice of rank $8$. 
 This implies 
\[ A_N\cong
\begin{cases}
  u(2)^{\oplus a }, \qquad\quad\qquad 0\le a \le 4,\ {\rm or} \\
  v(2)\oplus u(2)^{\oplus a }, \qquad 0\le a \le 3. 
\end{cases}
\]
 by the decomposition theorem of $2$-elementary finite quadratic forms (see 
 \cite[Proposition 1.8.1]{Nikulin1}).\\
 \quad Checking the signature of $(A_N,q_N)$, we see that the latter does 
 not occur. In all other cases, we find that $N$ has a unimodular overlattice
 of rank $8$, i.e., $E_8$. The index $[E_8 : N]$ is given by $2^a$.\\
$\claim$ {\itshape Let $N\subset E_8$ be an overlattice and assume 
 $N$ contains no $(-2)$ vectors.
 Then $[E_8:N]\ge 9$.}\\
$Proof.$ We take a basis of $E_8$ as in Figure \ref{fig2}. Consider 
 the elements $f_0:=0,\ f_j:=e_1+\cdots+e_j,\ 1\le j\le 7$ and $f_8:=
 2e_1+3e_2+4e_3+5e_4+6e_5+4e_6+2e_7+3e_8$. It is easy to see that any 
 difference $f_j-f_i,\ 0\le i<j\le 8$ has norm $-2$. This means that every $f_i$ 
 is a distinct element of the reminder class group $E_8/N$. $\square$
\\
 Thus we obtain $A_N=u(2)^{\oplus 4}.$ 
 Therefore $N$ has the maximal number of minimal generators 
 of the discriminant group. This implies that for any $n\in N$, $n/2\in N^*$.
 In particular $N(1/2)$ is a unimodular lattice.
 Thus we see that $N\cong E_8(2)$ and (1) is proved. \\
 \quad For (2), we use the notation introduced in Theorem \ref{lemb}, (1).
 Since $q_M\oplus q_N$ is nondegenerate
 $\#\Gamma=\#\Gamma_N=(\#A_M\cdot \#A_N/\#A_S)^{1/2}=2^7.$
 Therefore $\#\Gamma_N^{\perp}=\#A_N/\#\Gamma=2$ and we can put
 $\Gamma_N^{\perp}=\{0,z_N\}$.
 There are two cases where $q_N(z_N)=0,\ 1$.
 Thus there are at least two distinct primitive embedding of $M\subset S$.
 On the other hand, the canonical homomorphism $\sigma_M: O(M)\rightarrow O(q_M)$
 and $\sigma_N: O(N)\rightarrow O(q_N)$ are both surjective by Proposition \ref{useful}
 and \cite{B-P}. Thus the primitive embedding of $M\subset S$
 is classified by the invariant $q_N(z_N)$ as we see from the following argument.\\
 $\claim$ {\itshape Let $M_1$ and $M_2$ be two primitive sublattices of $S$ 
 satisfying the assumption. For each $k=1,2$, we use the same notation
 as above, indexed by $k$. If $q_{N_1}(z_{N_1})=q_{N_2}(z_{N_2})$, then there
 exists $\vp\in O(S)$ which transforms $M_1$ onto $M_2$.}\\  
 $Proof.$ By Witt's theorem on the finite quadratic forms, there exist $\psi_M$ and 
 $\psi_N$ fitting in the square inside.
\[\begin{CD}
 A_{M_1} \supset\ @. \Gamma_{M_1} @>\gamma_1>> \Gamma_{N_1} @. \ \subset A_{N_1}\\
@V\vp_M VV @V\psi_M VV @VV\psi_N V @VV\vp_N V \\
 A_{M_2} \supset\ @. \Gamma_{M_2} @>\gamma_2>> \Gamma_{N_2} @. \ \subset A_{N_2}
\end{CD}\]
\quad Again by Witt's theorem we can extend $\psi_M$ (resp. $\psi_N$) to
 $\vp_M$ (resp. $\vp_N$) in the diagram.
 By the surjectivity of $\sigma_M$ and $\sigma_N$ mentioned above, 
 these isomorphisms lift to an isomorphism between $M_1\oplus N_1$ and 
 $M_2\oplus N_2$ which preserves the overlattice $S$.
 This was the assertion. \qquad q.e.d.
\begin{rem}
 The same but geometric situation of the proposition
 is considered in \cite{M-N}. There (1) is proved by a geometric method.
\end{rem} 
 Now we regard $S$ as the N\'{e}ron-Severi lattice $NS$ of a generic
 Kummer surface of product type. To see the natural correspondence in 
 Theorem \ref{Kummer}, we have to associate a free involution with an 
 element of $A_{NS(X)}=A_{NS}$.
\begin{dfn}\label{invariant} 
 Let $M\cong U(2)\oplus E_8(2)\subset NS$ be a primitive sublattice.
 We use the same notation as in Proposition \ref{E82}.
 Then we define the {\it patching element}
 $v_M\in A_{NS}$ associated with the primitive sublattice $M$ to be the image of
 the element $(0,z_N)\in A_M\oplus A_N$ in $A_{NS}$ under the isomorphism of 
 Theorem \ref{lemb} (1).
\end{dfn}
 We note that $v_M\neq 0$. Also we note the equality $q_N(z_N)=q_{NS}(v_M)$. 
 Using $v_M$, we can describe the group $\sigma_{NS}(K)$ 
 in the Theorem \ref{bound}. Recall that $K\subset O(NS)$ is the stabilizer subgroup 
 of $M$. 
\begin{lem}
 $\sigma_{NS} (K)$ is equal to the stabilizer subgroup $G_{v_M}\subset O(q_{NS})$ of $v_M$.
\end{lem}
$Proof.$
 According to the isomorphism in Theorem \ref{lemb} (1), 
 we consider the group 
 \[H:=\{(\alpha_M,\alpha_N)\in O(q_M)\times O(q_N)|\alpha_M(\Gamma_M)=\Gamma_M
 ,\alpha_N(\Gamma_N)=\Gamma_N, \alpha_M \gamma =\gamma \alpha_N \}\]
 and decompose $\sigma_{NS}$ as
$\begin{CD}
 K @>p>> H @>q>> O(q_{NS}).
\end{CD}$
 Since an automorphism in $K$ preserves 
 $\Gamma_M$ and $\Gamma_N$,
 it is clear that $\sigma_{NS} (K)\subset G_{v_M}$.\\
 \quad We prove the converse, $\sigma_{NS} (K)\supset G_{v_M}$.
 First we note that $p$ is surjective since $\sigma_M$ and $\sigma_N$ are 
 both surjective as seen in the proof of Proposition \ref{E82}.
 Thus it is enough to see $\im q \supset G_{v_M}$.
 By Theorem \ref{lemb}, we have an isomorphism 
\begin{equation}
 q_M\cong (q_N\oplus q_{NS}|\Gamma'^{\perp})/\Gamma'  \label{discri}
\end{equation}
 where $\Gamma'$ is an isotropic subgroup of $q_N\oplus q_{NS}$ 
 which is a pushout of an isomorphism
 $\gamma': \Gamma'_N \rightarrow \Gamma'_{NS}$ between
 subgroups of $A_N$ and $A_{NS}$ respectively.
 Using the same notation of Proposition \ref{E82}, it is easy 
 to see that $\Gamma'=\{0, (z_N,v_M)\}$.\\
 \quad Suppose we are given an element $\beta\in G_{v_M}$.
 Then the automorphism 
\[ (\id_N,\beta)\in O(q_N)\times O(q_{NS})\]
 clearly preserves $\Gamma'$ and induces an element 
 $\alpha_M\in O(q_M)$ under the isomorphism (\ref{discri}).
 Since the element $(0,v_M)$ of the right-hand-side of (\ref{discri})
 corresponds to $z_M$, $\alpha_M$ preserves $z_M$.
 By construction $\alpha_M \gamma =\gamma \id_N$ holds 
 and therefore $(\alpha_M, \id_N)\in H$.
\qquad q.e.d.
\begin{prop}
 Let $M_1$ and $M_2$ be the fixed lattices of two free involutions 
 $i_{M_1},i_{M_2}$ on $X$. Then they give the isomorphic 
 quotients if and only if their patching elements coincide.
\end{prop}
$Proof.$ Suppose $X/i_{M_1}\cong X/i_{M_2}$. 
 There exists $\vp\in \Aut (X)$ such that
 $\vp(M_1)=M_2$. $\vp$ preserves the overlattice $NS$
 so that we have $\vp(z_{M_1})=z_{M_2}$ and $\sigma_{NS} (\vp)(v_{M_1})=v_{M_2}$.
 On the other hand, $\vp$ acts on $A_{NS}$ trivially by the
 assumption of the theorem. Thus we see $v_{M_1}=v_{M_2}$.\\
 \quad Conversely assume the patching elements $v_M\in A_{NS}$ concide.
 By Proposition \ref{E82}, the primitive embeddings $M_1$ and $M_2$ are 
 isomorphic and there exists $\vp\in O(NS)$ such that $\vp(M_1)=M_2$,
 namely $\vp i_{M_1} \vp^{-1}=i_{M_2}$.
 By assumption, $\sigma_{NS} (\vp)(v_M)=v_M$.
 We can assume $\vp\in O^{\uparrow}(NS)$ by replacing $\vp$ by $-\vp$ 
 if necessary. According to Proposition \ref{dec} (2),
 $\vp=w\psi$ where $w\in W_X,\psi \in O^+(NS)$. Then Lemma \ref{onFP} (2) 
 implies $\psi (M_1)=(M_2)$. On the other hand, since $w$ acts on 
 $A_{NS}$ trivially, $\sigma_{NS}(\psi)=\sigma_{NS}(\vp)$ and this element 
 fixes $v_M$. We apply Theorem \ref{bound}, Step 3 
 to $M=M_1, \psi_1=\psi, \psi_2=\id_{NS}$. Both $\psi_1$ and $\id$ stabilize
 $v_M$, therefore their images by $\sigma_{NS}$ have the same
 class in $O(q_{NS})/\sigma_{NS}(K)$ by the previous lemma.
 Thus the conclusion holds.  \qquad q.e.d.\\
\quad \\
 \quad Next we compute the patching elements of involutions of Lieberman and Kondo-Mukai.
 They involve parameters as mentioned in the beginning of this section and 
 we have to consider the dependence of patching elements on the parameters. 
 This is directly done. We take the following basis of $A_{NS}$.\\
 $e_1=(G_{11}+G_{13}+G_{31}+G_{33})/2$,\quad 
 $f_1=(G_{22}+G_{23}+G_{32}+G_{33})/2,$\\
 $e_2=(G_{21}+G_{23}+G_{31}+G_{33})/2,$\quad 
 $f_2=(G_{12}+G_{13}+G_{32}+G_{33})/2$.\\
 Then the result is as in Figures \ref{fig3},\ref{fig4}.
\quad \\
 \quad From these figures we can see which patching element we obtain 
 when we choose a parameter of a geometrically constructed free involution.
 In Figure \ref{fig4},
 we normalized the cases to only $i_4=j_4=4$.\\
 \quad Since the vectors in Figures \ref{fig3},\ref{fig4} run
 all over $A_{NS}-\{ 0 \}$, we obtain
\begin{prop}
 The two kinds of free involutions gives the all distinct Enriques quotients
 of $X$.
\end{prop} 
 The proof of Theorem \ref{Kummer} is completed.
\begin{rem}
 (1) The involution $\tau$ acts trivially on $NS(X)$. 
 So it induces a numerically trivial involution on the fifteen Enriques 
 quotients. The Kondo-Mukai case of this is the last and missing result of 
 \cite{M-N}, first found in \cite{Kondo}.\\
 (2) The number of Enriques quotients can be computed in other $2$-elementary cases,
 using the argument of this section.\\
 \quad When $NS(X)\cong U(2)\oplus E_8^{\oplus 2}$, the Barth-Peters case, 
 the number $B_0$ is equal to $1$.\\
 \quad When $NS(X)\cong U(2)\oplus E_8(2)$, then $X$ has only one
 Enriques quotient.\\
 Finally using the result of \cite{M-N}, 
 we see that in other $2$-elementary cases $X$ has no 
 Enriques quotients.\\
 (3) The generators of the whole automorphism group $\Aut (X)$ are 
 found in \cite{K-K}.\\
\end{rem}

\end{document}